\documentclass[a4paper]{amsart}

\usepackage{amsfonts}
\usepackage{amssymb}
\usepackage{amsfonts}
\usepackage{amssymb}
\usepackage{amsthm}
\usepackage{amsmath}

\usepackage{mathpple}

\theoremstyle{plain}

\newtheorem{teo}{Theorem}[section]

\newtheorem{prop}[teo]{Proposition}
\newtheorem{cor}[teo]{Corollary}

\theoremstyle{definition}

\newtheorem{conge}[teo]{Conjecture}

\theoremstyle{remark}

\newtheorem{rem}[teo]{Remark}
 
\newtheorem{es}[teo]{Example}

\numberwithin{equation}{section}

\def\R{{{\mathbb R}}}

\def\HH{{{\mathcal H}}}
\def\NN{{{\mathbb N}}}

\def\composed{\circ}
\def\comp{\composed}

\def\DG{{\mathcal G}_k}
\def\HHH{{\mathrm H}}

\def\qol{{q}}
\def\BBB{{\mathrm B}}

\begin{document}

\title[Some Properties of The Distance Function]
{Some Properties of the Distance Function and a
  Conjecture of De~Giorgi}
\author[Manolo Eminenti]{Manolo Eminenti}
\address[Manolo Eminenti]{Scuola Normale Superiore\\Piazza Cavalieri
  7\\56126 Pisa\\Italy}
\email[M. Eminenti]{manolo@linuz.sns.it}
\author[Carlo Mantegazza]{Carlo Mantegazza}
\address[Carlo Mantegazza]{Scuola Normale Superiore\\Piazza Cavalieri
7\\56126 Pisa\\Italy}
\email[C. Mantegazza]{mantegazza@sns.it}
\date{\today}

\begin{abstract} In the paper~\cite{degio5} 
Ennio~De~Giorgi conjectured 
that any compact $n$--dimensional regular 
submanifold $M$ of $\R^{n+m}$, moving by the gradient of the functional
$$
\int_M 1+\vert\nabla^{k}\eta^M\vert^2\,d\HH^n\,,
$$
where $\eta^M$ is the square of the distance function from the
submanifold $M$ and $\HH^n$ is the $n$--dimensional Hausdorff 
measure in $\R^{n+m}$, does not develop singularities in finite
time provided $k$ is large enough, depending on the dimension $n$.\\
We prove this conjecture by
means of the analysis of the geometric properties of the high
derivatives of the distance function from a submanifold of the
Euclidean space. In particular, we show 
some relations with the second fundamental form and its covariant
derivatives of independent interest.
\end{abstract}

\maketitle

\section{Introduction}

In the paper~\cite{degio5} (see also~\cite[Section~5]{degio6}) 
Ennio~De~Giorgi stated the following conjecture
(Congettura~2, Pag.~267):

\smallskip

{\em Any compact $n$--dimensional regular
  submanifold $M$ of $\R^{n+m}$ moving by the gradient of the
  functional
$$
\DG=\int_M
1+\vert\nabla^{k}\eta^M\vert^2\,d\HH^n\,,
$$
where $\eta^M$ is the square of the distance function from $M$ and
$\HH^n$ is the $n$--dimensional Hausdorff 
measure in $\R^{n+m}$, does not develop singularities in finite time
if $k>n+1$.}

\bigskip

We make some preliminary comments before proceeding with
the analysis. 

The regular submanifold $M$ can be described with an embedding 
$\varphi:M\to\R^{n+m}$ which induces a metric tensor $g$ on
$M$, by pulling back the standard scalar product of $\R^{n+m}$,
turning $(M,g)$ in a smooth Riemannian manifold isometrically embedded
in $\R^{n+m}$ via the map $\varphi$. Then, we let $\mu$ and $\nabla$ to
be respectively the canonical volume measure, which coincides with the
Hausdorff measure restricted to the image $\varphi(M)$, 
and the covariant differentiation operator on $(M,g)$.\\
Despite of the use of the same symbol, 
the iterated gradient which appears in the functional $\DG$ is not
such covariant differentiation but the standard $k$--order differential
in the canonical basis of $\R^{n+m}$.

When $k\geq3$ the variational gradient flow associated
to the functional $\DG$ is governed by a parabolic system of order
higher than 2, precisely of order $2k-2$ (see~\cite{ambman1}), 
hence, maximum principle and comparison theorems are not available. 
This means that initially embedded submanifolds can
possibly develop self--intersections during the flow.\\
Consequently, looking for more ``robust'' functionals $\DG$ in order to
deal also with {\em immersed}--only submanifolds
$\varphi:M\to\R^{n+m}$, two problems naturally arise, one is the
difference (because of the possible multiplicities) between the
Hausdorff measure on the image $\varphi(M)$ and the canonical volume
measure $\mu$ on $(M,g)$ which can be overcome substituting $\HH^n$ 
with this latter, the second point is the non smoothness of the distance function
near the points of self--intersection.\\
If a compact smooth submanifold is embedded, the square of the 
distance function from $M$ turns out to be locally
smooth so the computation of the derivatives gives no trouble. 
This is no more true when $M$ is only immersed, hence, in such 
situation, at any $p\in M$ we will consider the derivatives of the distance 
function from an {\em embedded} image $\varphi(U_p)$ of a local
neighborhood $U_p\subset M$ of the point $p$.\\
Thus, we redefine the functionals as follows
\begin{equation}\label{funct}
\DG(\varphi)=\int_M
1+\vert\nabla^{k}\eta^M\vert^2\,d\mu\,,
\end{equation}
on the space of smooth immersions $\varphi:M\to\R^{n+m}$ of a compact
$n$--dimensional manifold $M$, where $\eta^M$ is the square of the
distance function from $\varphi(M)$ (keeping in mind the previous
discussion at the points of self intersections) 
and $\mu$ is the canonical volume measure associated to
the Riemannian manifold $(M,g)$, with
$g=\varphi^*\langle\,,\,\rangle_{\R^{n+m}}$.\\
Then, De~Giorgi's conjecture can be restated as follows,

\begin{conge}\label{dgcong} Any initial $n$--dimensional smooth
  submanifold $\varphi:M\to\R^{n+m}$ evolves by the gradient of the
  functional $\DG$ without developing singularities in finite time, if $k>n+1$.
\end{conge}

We will see that actually the weaker hypothesis
$k>[n/2]+2$, where $[n/2]$ denotes the integer part of $n/2$, is
sufficient.

\smallskip

In order to show such conjecture, we work out some properties, of
independent interest, about the high derivatives of the square of the distance
function from a submanifold, in particular their relation with the
covariant derivatives of the second fundamental form. This is the 
goal of the first part of the paper which can be seen as a
continuation of the analysis carried out in~\cite{ambman1}.\\
Then, in Section~\ref{prova} we show that these properties imply 
a priori estimates on the Sobolev constants of the evolving
manifolds, which allow us to follow the method used 
in~\cite{mant5} to prove the regularity of the flow
associated to the functionals
$$
{\mathcal F}_k(\varphi)=\int_M1+\vert\nabla^k\nu\vert^2\,d\mu\,,
$$
where $\nu$ is the normal vector field of a hypersurface $M$ in the
Euclidean space.

We conclude the paper discussing the subsequent open problem of
De~Giorgi, again stated in~\cite{degio5} (see also~\cite{degio6}),
about the singular approximation of the motion
by mean curvature with these smooth higher order flows.

\section{The Squared Distance Function from a Manifold}

We denote with $e_1,\dots,e_{n+m}$ the canonical basis of
$\R^{n+m}$ and with $\langle\,,\,\rangle$ its standard scalar product.\\
We let $M\subset\R^{n+m}$ be a smooth, compact,
$n$--dimensional, regular submanifold without boundary, then 
$T_xM$ and $N_xM\subset\R^{n+m}$ are, respectively, the tangent space
and the normal space to $M$ at $x\in M\subset\R^{n+m}$.

The {\em distance function} $d^M(x)$ and the
{\em squared distance function} 
$\eta^M(x)$ from $M$ are simply given by
$$
d^M(x)=\inf_{y\in M}\vert x-y\vert\qquad\text{ and }\qquad 
\eta^M(x)=[d^M(x)]^2
$$
for any $x\in\R^{n+m}$ (we will drop the superscript $M$ when no ambiguity
is possible). In this section we recall some facts from~\cite{ambman1}
about the distance function and we establish some new relations
between the high derivatives of $\eta^M$ and the second fundamental
form of $M$.

Since $M$ is smooth, embedded and compact, there exists an open
neighborhood $\Omega\subset\R^{n+m}$ of $M$ such that 
$d^M$ is smooth in $\Omega\setminus M$ and $\eta^M$ is smooth in
all $\Omega$.\\
Clearly, $\eta^M$ and $\nabla\eta^M(x)=0$ at every point of $x\in M$,
moreover, for every $x\in\Omega$ we have that $x-\nabla\eta^M(x)/2$ is
{\em the unique} point in $M$ of minimum distance from $x$ (the {\em
  projection} of $x$ on $M$), that we denote with $\pi^M(x)$.\\
Another nice property of the squared distance is that, for every 
$x\in M$ the Hessian matrix $\nabla^2\eta^M(x)$ is twice the matrix of
orthogonal projection onto the normal space $N_xM$. We will denote
respectively with $X^M$ and $X^\perp$ the projections of a vector $X$
on the tangent and normal space of $M$.

Let $x\in M$ and $X,\,Y\in T_x M$, the vector valued 
{\em second fundamental form} of $M$ at the point $x$ is given by 
$$
\BBB(X,Y)=\Bigl(\frac{\partial Y(x)}{\partial X}\Bigr)^\perp
$$
where we extended locally the two vectors $X$, $Y$ to 
tangent vector fields on $M$ (the derivative is
well defined since $X$ is a tangent vector at $x$).\\
If $\{\nu_\alpha\}_{\alpha=1,\dots, m}$ is a local basis of the normal
  bundle we have clearly
$$
\BBB(X,Y)=-\sum_{\alpha=1}^m \Bigl\langle\frac{\partial
  \nu_\alpha(x)}{\partial X},Y\Bigr\rangle \nu_\alpha\,.
$$
We will see $\BBB$ as a bilinear map from $T_xM\times T_xM$ to
$\R^{n+m}$, hence, as a family of $n+m$ 
bilinear forms $\BBB^k=\langle\BBB,e_k\rangle: T_xM\times
T_xM\to\R^{n+m}$. Moreover, we consider $\BBB$ acting also on vectors of
$\R^{n+m}$, not necessarily tangent, by setting 
$\BBB(V,W)=\BBB(V^M,W^M)\in N_xM\subset \R^{n+m}$ for every pair 
$V,\, W\in\R^{n+m}$. With such definition,
$\BBB^k_{ij}=\langle\BBB(e_i,e_j),e_k\rangle$.\\
It is well known that ${\BBB}$ is a symmetric bilinear form and its
trace is the {\em mean curvature} of components
${\HHH^k}=\sum_j\BBB^k_{jj}$.

We introduce now the function
$$
A^M(x)=\frac{\vert x\vert^2-[d^M(x)]^2}{2}\,,
$$
smooth as $\eta^M$ in the neighborhood $\Omega$ of $M$, and we set
$$
A^M_{i_1 \dots i_k}(x)=\frac{\partial^kA^M(x)}{\partial x_{i_1}\dots\partial
  x_{i_k}}
$$
for the derivatives of $A^M$ at every point $x\in\Omega$.\\
The  following Proposition (see~\cite{ambman1} for the proof) shows the first 
connection between the second fundamental form and the function $A^M$
(or equivalently, the squared distance function).

\begin{prop} The following relations hold,
\begin{itemize}
\item For any $x\in\Omega$, the point $\nabla A^M(x)$ is the projection
point $\pi^M(x)$.
\item If $x\in M$, then $\nabla^2A^M(x)$ is the matrix of orthogonal projection on $T_xM$.
\item For every $x\in M$, 
\begin{eqnarray*}
\BBB^k_{ij}&=&A^M_{ijs}(\delta_{ks}- A^M_{ks})\\
A^M_{ijk}&=&\BBB^k_{ij}+\BBB^i_{jk}+\BBB^j_{ki}\\
\HHH^k&=&\sum_i A^M_{iik}
\end{eqnarray*}
\end{itemize}
\end{prop}

We define now the $k$--derivative tensor $A^{k}(x)$ working on the $k$--uple
of vectors $X_i\in\R^{n+m}$, where $X_i=X_i^je_j$, as follows
$$
A^{k}(x)(X_1, \dots, X_k)=A^M_{i_1 \dots i_k}(x)X_1^{i_1}\dots
X_k^{i_k}\,,
$$
notice that the tensors $A^k$ are symmetric.\\
By sake of simplicity, we dropped the superscript $M$ on $A^k$, by the
same reason, we will also avoid to indicate the point $x\in M$ in the
sequel.

Our goal is to express $A^{k}$ in terms of covariant derivatives 
of $\BBB$.

\begin{prop}\label{teoremauno}
For every $k\geq2$ and for every $s\in\{1,\dots k\}$ there exists a
family $p^{k,s}_{j_1 \dots j_{k-s}}$ of symmetric polynomial tensors of type
$(s,0)$ on $M$, where ${j_1, \dots,j_{k-s}\in\{1,\dots,n+m\}}$, 
which are contractions of the second fundamental form $\BBB$
and its covariant derivatives with the metric tensor $g$, such that
$$
A^k(X_1,\dots, X_s,N_1,\dots, N_{k-s})= p^{k,s}_{j_1 \dots j_{k-s}}(X_1,\dots,
X_s)N_1^{j_1} \dots N_{k-s}^{j_{k-s}}
$$
for every $s$--uple of tangent vectors $X_h$ and $(k-s)$--uple of
normal vectors $N_h=N^{j_h}_he_h\in\R^{n+m}$.\\
Moreover, the tensors $p^{k,s}_{j_1 \dots j_{k-s}}$ are invariant by
exchange of the $j$--indices and the maximum order of differentiation  of
$\BBB$ which appears in every $p^{k,s}_{j_1 \dots j_{k-s}}$ is at most
$k-3$.
Considering the tangent plane at any point $x\in M$ also as a subset of
$\R^{n+m}$, the polynomials tensors $p^{k,s}_{j_1 \dots j_{k-s}}$ 
are expressed in the coordinate basis of the Euclidean space as follows 
$$
p^{k,s}_{j_1 \dots j_{k-s}}(X_1,\dots,
X_s)N_1^{j_1} \dots N_{k-s}^{j_{k-s}}=
p^{k,s}_{j_1 \dots j_{k-s}, i_1 \dots i_s}X_1^{i_1} \dots X_s^{i_s}
  N_1^{j_1} \dots N_{k-s}^{j_{k-s}}\,.
$$
Then, a family of tensors satisfying the above properties can be
defined recursively according to the following formulas
\begin{align}
p^{2,0}_{j_1j_2}=\,&p^{2,1}_{j_1,i_1}=0\,,\qquad p^{2,2}_{i_1i_2}=\delta_{i_1i_2}\label{p2line}\\
p^{k,0}_{j_1\dots j_k}=\,&p^{k,1}_{j_1\dots
  j_{k-1},i_1}=0\qquad\qquad\qquad\qquad\qquad&\text{ for every
  $k\geq 2$}\label{k01line}\\
p^{k+1,s}_{j_1 \dots j_{k-s+1}, i_0 i_1 \dots i_{s-1}}
=\,&(\nabla p^{k,s-1}_{j_1 \dots j_{k-s+1}})_{i_0 i_1 \dots
  i_{s-1}}&\text{{ if $s<k+1$}}\label{induc}\\
&-\sum_{h=1}^{k-s+1} p^{k,s-1}_{j_1 \dots j_{h-1}rj_{h+1}\dots
  j_{k-s+1},i_1 \dots i_{s-1}}\BBB_{ri_0}^{j_h}\nonumber\\
&-\sum_{h=1}^{s-1}p^{k,s-2}_{rj_1 \dots j_{k-s+1},i_1 \dots i_{h-1} i_{h+1} \dots
  i_{s-1}}\BBB^{r}_{i_0i_h}\nonumber\\
&+\sum_{h=1}^{k-s+1} p^{k,s}_{j_1 \dots j_{h-1}j_{h+1}\dots
  j_{k-s+1},i_1 \dots i_{s-1}r} \BBB_{ri_0}^{j_h}\nonumber\\
p^{k+1,k+1}_{i_0 i_1 \dots i_{k+1}}
=\,&-\sum_{h=1}^{k}p^{k,k-1}_{r,i_1 \dots i_{h-1} i_{h+1} \dots
  i_{k}}\BBB^{r}_{i_0i_h}\,.\label{lastline}
\end{align}
\end{prop}

\begin{proof}
If $k=2$ we have immediately
$$ 
A^{2}(N_1,N_2)=0, \qquad A^{2}(X_1,N_1)=0, \qquad
A^{2}(X_1,X_2)=X^i_1X^i_2=\delta_{i_1i_2}X_1^{i_1}X_2^{i_2}
$$
since $X_1$ and $X_2$ are tangent and $A^{2}$ is the projection on
the tangent space. Hence, formula~\eqref{p2line} follows.\\
We argue now by induction on $k\geq2$. When $s=0$ the value 
$A^{k}(N_1,\dots,N_k)(x)$ depends only on the function $A^M$ on the
$m$--dimensional normal subspace to $M$ at $x$, and on this subspace
$A^M$ is identically zero, hence the first equality in~\eqref{k01line}
is proved.\\
Suppose now that $s\in\{1,\dots, k+1\}$, we extend the vectors $X_h\in
T_xM$ and $N_h\in N_xM$ to a family of local vector fields,
respectively tangent and normal to $M$, then
\begin{align*}
A^{k+1}(X_0,X_1, \dots, &\,X_{s-1},N_1,\dots, N_{k-s+1})=
\frac{\partial\,}{\partial X_0}\Bigl(A^{k}(X_1,
\dots,X_{s-1},N_1,\dots, N_{k-s+1})\Bigr)\\
&\,-\sum_{h=1}^{s-1}A^{k}\Bigl(X_1,\dots X_{h-1},\frac{\partial
  X_h}{\partial X_0},X_{h+1},\dots,X_{s-1},N_1,\dots, N_{k-s+1}\Bigr)\\
&\,-\sum_{h=1}^{k-s+1}A^{k}\Bigl(X_1,\dots,X_{s-1},N_1,\dots, 
\frac{\partial N_h}{\partial X_0}, \dots, N_{k-s+1}\Bigr)
\end{align*}
where the last line is not present in the special case $s=k+1$ and the
second line is not present if $s=1$. In this last case, we have
$$
A^{k+1}(X_0,N_1,\dots, N_{k})=
\frac{\partial\,}{\partial X_0}\Bigl(A^{k}(N_1,\dots, N_{k})\Bigr)
-\sum_{h=1}^{k}A^{k}\Bigl(N_1,\dots, 
\frac{\partial N_h}{\partial X_0}, \dots, N_{k}\Bigr)=0
$$
since the first term of the right member is zero by the first equality
in~\eqref{k01line} and, after decomposing $\frac{\partial
  N_h}{\partial X_0}$ in tangent and normal part, the tangent term is
zero by induction and the normal term is zero for~\eqref{k01line}
again. This shows the second equality in~\eqref{k01line}.\\
So we suppose $1<s<k+1$, by the inductive hypothesis, 
$$
A^{k}(X_1,\dots,X_{s-1}, N_1,\dots, N_{k-s+1})=
p^{k,s-1}_{j_1 \dots j_{k-s+1}}(X_1,\dots,
X_{s-1})N_1^{j_1} \dots N_{k-s+1}^{j_{k-s+1}}
$$
thus, differentiating along $X_0$, 
which is a tangent field, we obtain
\begin{align*}
A^{k+1}(X_0,&X_1, \dots, X_{s-1},N_1,\dots, N_{k-s+1})\\
=&\,\frac{\partial\,}{\partial X_0}
\Bigl(p^{k,s-1}_{j_1 \dots j_{k-s+1}}(X_1,\dots,
X_{s-1})N_1^{j_1} \dots N_{k-s+1}^{j_{k-s+1}}\Bigr)\\
&\,-\sum_{h=1}^{s-1}A^{k}\Bigl(X_1,\dots, \Bigl(\frac{\partial
  X_h}{\partial X_0}\Bigr)^M,\dots, X_{s-1}, N_1,\dots, N_{k-s+1}\Bigr)\\
&\,-\sum_{h=1}^{s-1}A^{k}\Bigl(X_1,\dots, \Bigl(\frac{\partial
  X_h}{\partial X_0}\Bigr)^\perp, \dots, X_{s-1},N_1,\dots, N_{k-s+1}\Bigr)\\
&\,-\sum_{h=1}^{k-s+1}A^{k}\Bigl(X_1,\dots,X_{s-1},N_1,\dots, 
\Bigl(\frac{\partial N_h}{\partial X_0}\Bigr)^M, \dots, N_{k-s+1}\Bigr)\\
&\,-\sum_{h=1}^{k-s+1}A^{k}\Bigl(X_1,\dots,X_{s-1},N_1,\dots, 
\Bigl(\frac{\partial N_h}{\partial X_0}\Bigr)^\perp, \dots,
N_{k-s+1}\Bigr)\,.
\end{align*}
We use now the symmetry of $A^k$ and 
we substitute recursively $p^{k,s}$, $p^{k,s-1}$ and
$p^{k,s-2}$ to $A^{k}$, according to the number of tangent vectors
inside $A^k$,
\begin{align*}
A^{k+1}(X_0,&X_1, \dots, X_{s-1},N_1,\dots, N_{k-s+1})\\
=&\,\frac{\partial\,}{\partial X_0}
\Bigl(p^{k,s-1}_{j_1 \dots j_{k-s+1}}(X_1,\dots,
X_{s-1})\Bigr) N_1^{j_1} \dots N_{k-s+1}^{j_{k-s+1}}\\
&\,+\sum_{h=1}^{k-s+1} p^{k,s-1}_{j_1 \dots j_{k-s+1}}(X_1,\dots,
X_{s-1}) N_1^{j_1} \dots \frac{\partial N_h^{j_h}}{\partial
  X_0} \dots N_{k-s+1}^{j_{k-s+1}}\\
&\,-\sum_{h=1}^{s-1} p^{k,s-1}_{j_1 \dots j_{k-s+1}}(X_1,\dots, 
\nabla_{X_0}X_h,\dots,X_{s-1}) N_1^{j_1} \dots N_{k-s+1}^{j_{k-s+1}}\\
&\,-\sum_{h=1}^{s-1}p^{k,s-2}_{rj_1 \dots j_{k-s+1}}(X_1,\dots,
X_{h-1},X_{h+1},\dots, X_{s-1}) \Bigl[\Bigl(\frac{\partial
  X_h}{\partial X_0}\Bigr)^\perp\Bigr]^{r}
N_1^{j_1}\dots N_{k-s+1}^{j_{k-s+1}}\\
&\,-\sum_{h=1}^{k-s+1} p^{k,s}_{j_1 \dots j_{h-1}j_{h+1}\dots
  j_{k-s+1}}\Bigl(X_1,\dots,X_{s-1},  \Bigl(\frac{\partial N_h}{\partial
  X_0}\Bigr)^M\Bigr)N_1^{j_1} \dots
N_{h-1}^{j_{h-1}}N_{h+1}^{j_{h+1}}\dots N_{k-s+1}^{j_{k-s+1}}\\
&\,-\sum_{h=1}^{k-s+1}p^{k,s-1}_{j_1 \dots j_{k-s+1}}(X_1,\dots,X_{s-1}) N_1^{j_1}\dots 
\Bigl[\Bigl(\frac{\partial N_h}{\partial X_0}\Bigr)^\perp\Bigr]^{j_h} \dots
N_{k-s+1}^{j_{k-s+1}}\,.
\end{align*}
Adding the first and the third line we get the covariant derivative of
$p^{k,s-1}_{j_1 \dots j_{k-s+1}}$, adding the second and the last line we
get
\begin{align*}
A^{k+1}(X_0,&X_1, \dots, X_{s-1},N_1,\dots,N_{k-s+1})\\
=&\,\nabla p^{k,s-1}_{j_1 \dots j_{k-s+1}}(X_0, X_1,\dots,
X_{s-1}) N_1^{j_1} \dots N_{k-s+1}^{j_{k-s+1}}\\
&\,+\sum_{h=1}^{k-s+1} p^{k,s-1}_{j_1 \dots j_{k-s+1}}(X_1,\dots,
X_{s-1}) N_1^{j_1} \dots \Bigl[\Bigl(\frac{\partial N_h}{\partial
  X_0}\Bigr)^M\Bigr]^{j_h}\dots N_{k-s+1}^{j_{k-s+1}}\\
&\,-\sum_{h=1}^{s-1}p^{k,s-2}_{rj_1 \dots j_{k-s+1}}(X_1,\dots,
X_{h-1},X_{h+1},\dots,X_{s-1}) \Bigl[\Bigl(\frac{\partial
  X_h}{\partial X_0}\Bigr)^\perp\Bigr]^{r} 
N_1^{j_1}\dots N_{k-s+1}^{j_{k-s+1}}\\
&\,-\sum_{h=1}^{k-s+1} p^{k,s}_{j_1 \dots j_{h-1}j_{h+1}\dots
  j_{k-s+1}}\Bigl(X_1,\dots,X_{s-1},  \Bigl(\frac{\partial N_h}{\partial
  X_0}\Bigr)^M\Bigr)N_1^{j_1} \dots
N_{h-1}^{j_{h-1}}N_{h+1}^{j_{h+1}}\dots N_{k-s+1}^{j_{k-s+1}}\,.
\end{align*}
Taking now into account  that 
$$
\Bigl[\Bigl(\frac{\partial N_h}{\partial X_0}\Bigr)^M\Bigr]^{r}=
-\BBB_{ri_0}^{j_h}X_0^{i_0}N_h^{j_h}
\text{ \,and \,} \Bigl[\Bigl(\frac{\partial   X_h}{\partial
  X_0}\Bigr)^\perp\Bigr]^r
= \BBB^r_{i_0i_h}X_0^{i_0}X_h^{i_h}
$$
we get
\begin{align*}
A^{k+1}(X_0,&X_1, \dots, X_{s-1},N_1,\dots,N_{k-s+1})\\
=&\, \nabla p^{k,s-1}_{j_1 \dots j_{k-s+1}}(X_0, X_1,\dots,
X_{s-1}) N_1^{j_1} \dots N_{k-s+1}^{j_{k-s+1}}\\
&\,-\sum_{h=1}^{k-s+1} p^{k,s-1}_{j_1 \dots j_{h-1}rj_{h+1}\dots j_{k-s+1}}(X_1,\dots,
X_{s-1})\BBB_{ri_0}^{j_h}X_0^{i_0}N_1^{j_1}\dots N_{k-s+1}^{j_{k-s+1}}\\
&\,-\sum_{h=1}^{s-1}p^{k,s-2}_{rj_1 \dots j_{k-s+1}}(X_1,\dots,
X_{h-1},X_{h+1},\dots,X_{s-1})\BBB^{r}_{i_0i_h}X_0^{i_0} X_h^{i_h}
N_1^{j_1}\dots N_{k-s+1}^{j_{k-s+1}}\\
&\,+\sum_{h=1}^{k-s+1} p^{k,s}_{j_1 \dots j_{h-1}j_{h+1}\dots
  j_{k-s+1}}(X_1,\dots,X_{s-1},  \BBB_{ri_0}^{j_h}
X_0^{i_0}e_r)N_1^{j_1}\dots N_{k-s+1}^{j_{k-s+1}}\,.
\end{align*}
Then, expressing the tensors in coordinates, we have
\begin{align*}
A^{k+1}(X_0,&X_1, \dots, X_{s-1},N_1,\dots, N_{k-s+1})\\
=&\,(\nabla p^{k,s-1}_{j_1 \dots j_{k-s+1}})_{i_0 i_1 \dots i_{s-1}}
X_0^{i_0}  X_1^{i_1}\dots X_{s-1}^{i_{s-1}} N_1^{j_1} \dots N_{k-s+1}^{j_{k-s+1}}\\
&\,-\sum_{h=1}^{k-s+1} p^{k,s-1}_{j_1 \dots j_{h-1}rj_{h+1}\dots
  j_{k-s+1},i_1 \dots i_{s-1}}\BBB_{ri_0}^{j_h}
X_0^{i_0}\dots X_{s-1}^{i_{s-1}}  N_1^{j_1}\dots N_{k-s+1}^{j_{k-s+1}}\\
&\,-\sum_{h=1}^{s-1}p^{k,s-2}_{rj_1 \dots j_{k-s+1},i_1 \dots i_{h-1} i_{h+1} \dots
  i_{s-1}}\BBB^{r}_{i_0i_h} X_0^{i_0} \dots X_{s-1}^{i_{s-1}} 
N_1^{j_1}\dots N_{k-s+1}^{j_{k-s+1}}\\
&\,+\sum_{h=1}^{k-s+1} p^{k,s}_{j_1 \dots j_{h-1}j_{h+1}\dots
  j_{k-s+1},i_1 \dots i_{s-1}r}\BBB_{ri_0}^{j_h} 
X_0^{i_0}\dots X_{s-1}^{i_{s-1}} N_1^{j_1} \dots N_{k-s+1}^{j_{k-s+1}}\,,
\end{align*}
which is formula~\eqref{induc}.\\
In the special case $s=k+1$, formula~\eqref{lastline}, we just have to
repeat the computations dropping all the lines containing sums like 
$\sum_{h=1}^{k-s+1}... $, which are not present.\\
Finally, assuming inductively that the polynomial tensors $p^{k,s}$, 
$p^{k,s-1}$ and $p^{k,s-2}$ are symmetric in the $j$--indices 
and contain covariant derivatives of $\BBB$ only up to the order
$k-3$,  also the claims about the symmetry and the order of the
derivatives of $\BBB$ follow.\\
\end{proof}

\begin{es}\label{exmp}
We compute some $p^{k,s}$ as a consequence of this proposition.

\begin{enumerate}
\item When $k=2$ we saw that
$$
p^{2,0}_{j_1j_2}=0,\qquad p^{2,1}_{j_1}=0,\qquad p^{2,2}=g\,.
$$
\item When $k=3$ we have
\begin{align*}
&\,p^{3,0}_{j_1j_2j_3}=0,\qquad p^{3,1}_{j_1j_2}=0,\\
&\,p^{3,2}_{j_1,i_1i_2}=p^{2,2}_{i_2r}\BBB_{ri_1}^{j_1}
=\BBB_{i_1i_2}^{j_1}\\
&\,p^{3,3}_{i_1i_2i_3}=(\nabla p^{2,2})_{i_1i_2i_3} +
p^{2,1}_{r,i_2}\BBB_{i_1i_3}^r+
p^{2,1}_{r,i_3}\BBB_{i_1i_2}^r=0
\end{align*}
that is,
$$
p^{3,2}_{j_1}=\BBB^{j_1}\text{ and }\, p^{3,3}=0\,.
$$
\item When $k=4$ we have
\begin{align*}
p^{4,0}_{j_1j_2j_3j_4}=&\,0,\qquad p^{4,1}_{j_1j_2j_3}=0,\\
p^{4,2}_{j_1j_2,i_1i_2}=&\, 
p^{3,2}_{j_1,i_1r}\BBB_{ri_2}^{j_2}
+p^{3,2}_{j_2,i_1r}\BBB_{ri_1}^{j_1}=
\BBB_{i_1r}^{j_1}\BBB_{ri_2}^{j_2}
+\BBB_{i_2r}^{j_2}\BBB_{ri_1}^{j_1}\,,\\
p^{4,3}_{j_1,i_1i_2i_3}= &\,(\nabla p^{3,2}_{j_1})_{i_1i_2i_3}
+p^{3,2}_{r,i_2i_3}\BBB_{ri_1}^{j_1}=
(\nabla p^{3,2}_{j_1})_{i_1i_2i_3}
+\BBB^{r}_{i_2i_3}\BBB_{ri_1}^{j_1}
=(\nabla \BBB^{j_1})_{i_1i_2i_3}\\
&\,\text{since we contracted a normal vector with a tangent one,}\\
p^{4,4}_{i_1i_2i_3i_4}=&\,
-p^{3,2}_r{i_3i_4}\BBB_{i_1i_2}^r
-p^{3,2}_r{i_2i_4}\BBB_{i_1i_3}^r
-p^{3,2}_r{i_2i_3}\BBB_{i_1i_4}^r
=-\BBB^r_{i_3i_4}\BBB_{i_1i_2}^r
-\BBB^r_{i_2i_4}\BBB_{i_1i_3}^r
-\BBB^r_{i_2i_3}\BBB_{i_1i_4}^r\,.
\end{align*}
\end{enumerate}
\end{es}

Proposition~\ref{teoremauno} allows us to write $A^k$ in terms of the
tensors $p^{k,s}$ and the projections on the tangent and normal spaces (hence
contracting with the scalar product of $\R^{n+m}$), so we get the
following corollary.

\begin{cor}\label{corouno}
For every $k\geq3$ the symmetric tensor $A^k$ can be expressed as a 
polynomial tensor in $\BBB$ and its covariant derivatives, contracted
with the scalar product of $\R^{n+m}$.\\
The maximum order of differentiation of $\BBB$ which appears
in $A^k$ is $k-3$. Precisely, the only tensors among the $p^{k,s}$ 
containing such highest derivative are $p^{k,k-1}_{j_1}$, given by
$$
p^{k,k-1}_{j_1}=\nabla^{k-3}\BBB^{j_1} + \,{\mathrm {LOT}}\,.
$$
where we denoted with ${\mathrm {LOT}}$ ({\em lower order terms}) 
a polynomial term containing only derivatives of $\BBB$ up to the order k-4.
\end{cor}
\begin{proof}
Looking at the tensors with the derivative of $\BBB$ of maximum
order among the $p^{k,s}_{j_1\dots j_{k-s}}$, 
by formula~\eqref{induc} and the fact that the only non zero 
polynomials $p^{3,s}_{j_1\dots j_{3-s},i_1\dots i_s}$ are 
$p^{3,2}_{j_1,i_1i_2}=\BBB^{j_1}_{i_1i_2}$ 
(see Example~\ref{exmp}), it is clear that they come from the
derivative $\nabla p^{k-1,k-2}_{j_1}$. Iterating the argument, 
the leading term in $p^{k,k-1}_{j_1}$ is given by
$\nabla^{k-3}p^{3,2}_{j_1}=\nabla^{k-3}\BBB^{j_1}$.
\end{proof}
\begin{rem}
We can see in Example~\ref{exmp} that 
when $k=3$ and $4$, the lower order term which appears 
above is zero, actually, when $k\geq 5$ this is no more true.
\end{rem}

The decomposition of $A^k$ in its tangent and normal components is
very useful in studying the norm of $A^k$, which is the main quantity
we are interested in.

Fixing at a point $x\in M$ an orthonormal basis $\{e_1,\dots,
e_{n+m}\}$ of $\R^{n+m}$ such that $\{e_1,\dots, e_n\}$ is a basis of
$T_xM$, we have obviously
\begin{align*}
\vert A^k\vert^2&\,=\sum_{1\leq i_1,\dots,i_k\leq
  n+m}[A^k(e_{i_1},\dots,e_{i_k})]^2\\
&\,\geq \sum_{{\genfrac{}{}{0pt}{}{{1\leq i_1, i_2\leq n}}
{{n< i_3 ,\dots, i_k \leq n+m}}}}
[A^k(e_{i_1},e_{i_2},e_{i_3},\dots,e_{i_k})]^2\\
&\, \geq \sum_{n<j\leq n+m}\sum_{1\leq i_1, i_2\leq n}
[A^k(e_{i_1}, e_{i_2},e_j,\dots,e_j)]^2\\
&\, =\sum_{n<j\leq n+m}\sum_{1\leq i_1, i_2\leq n}
[p^{k,2}_{j\dots j,i_1i_2}]^2\,,
\end{align*}
that is,
$$
\vert A^k\vert^2\geq \sum_{n< j\leq n+m}
\vert p^{k,2}_{j\dots j}\vert^2\,.
$$
We analyse this last term by means of formula~\eqref{induc}. We
have $p^{2,2}=g$ and for every $k\geq 2$, 
$$
p^{k+1,2}_{j\dots j, i_0 i_1}=
\sum_{h=1}^{k-1} p^{k,2}_{j\dots j,i_1r}\BBB_{ri_0}^{j}
=(k-1)\,p^{k,2}_{j\dots j,i_1r}\BBB_{ri_0}^{j}\,.
$$
Then, by induction, it is easy to see that
$$
p^{k,2}_{j\dots j, i_0 i_1}=(k-2)!\,\BBB_{i_0r_1}^{j}\BBB_{r_1r_2}^j\dots\BBB_{r_{k-3}i_1}^j
$$
hence, as the bilinear form $\BBB^j$ is symmetric, denoting with
$\lambda^j_s$ its eigenvalues at the point $x\in M$, we conclude
$$
\vert p^{k,2}_{j\dots j}\vert^2=[(k-2)!]^2\sum_{s=1}^n (\lambda^j_s)^{2(k-2)}
\geq \widetilde{C}_k\vert\BBB^j\vert^{2k-4}\,.
$$
Coming back to our estimate,
$$
\vert A^k\vert^2\geq \widetilde{C}_k\sum_{n< j\leq n+m}
\vert\BBB^j\vert^{2k-4}\geq 
C_k\Bigl(\sum_{n< j\leq n+m}\vert\BBB^j\vert^2\Bigr)^{k-2}
=C_k\vert\BBB\vert^{2k-4}\,.
$$

\begin{prop}\label{akteo} The following estimate holds,
\begin{equation*}
\vert A^k\vert^2\geq C_k\vert\BBB\vert^{2k-4}
\end{equation*}
where $C_k$ is an universal constant depending only on $k$, $n$ and
$m$.
\end{prop}

\section{Flows by Geometric Functionals}\label{prova}

The very first step in proving De~Giorgi's Conjecture~\ref{dgcong} is
to see that any initial submanifold actually moves smoothly by the
gradient of the functional $\DG$ in~\eqref{funct}, at least for
some small time.

In the paper~\cite{ambman1}, Theorem~4.5
and Theorem~5.9, it has been shown that the 
{\em first variation} of the functional
$\DG$ is given by
\begin{equation*}
E_{\DG}=\Bigl(-\HHH + 2k(-1)^{ k -1}
\overset{\text{$( k -2)$--times}}{\overbrace{\Delta^M\comp
\Delta^M\comp\,\dots\,\comp\Delta^M}}\,\,\HHH\Bigr)^\perp
+h^j(A^M)e_j^\perp
\end{equation*}
where the functions $h^j(A^M)$ are polynomials in the derivatives of
$A^M$ up to the order $2k -2$.\\
By means of Corollary~\ref{corouno} we can express the terms $h^j(A^M)$
as polynomials $q^j(\BBB)$ obtained contracting 
$\BBB$ and its covariant derivatives up to the order $2k-5$ with the
scalar product of $\R^{n+m}$.\\
To get a solution of the geometric evolution problem for the initial
submanifold $\psi:M\to\R^{n+m}$, we look then 
for a smooth function $\varphi:M\times[0,T)\to\R^{n+m}$ such that
\begin{enumerate}
\item the map $\varphi_t=\varphi(\cdot,t):M\to\R^{n+m}$ is an
  immersion for every $t\in[0,T)$;
\item $\varphi_0(p)=\varphi(p,0)=\psi(p)$ for every $p\in M$;
\item the following parabolic system is satisfied
\begin{equation*}
\frac{\partial\varphi}{\partial t}=\HHH+ 2k(-1)^{k}
\Bigl(\overset{\text{$( k -2)$--times}}{\overbrace{\Delta^{M_t}\comp
\Delta^{M_t}\comp\,\dots\,\comp\Delta^{M_t}}}\,\,\HHH\Bigr)^\perp
+q^j(\BBB)e_j^\perp\,.
\end{equation*}
\end{enumerate}
Here we denoted with $\Delta^{M_t}$ the Laplacian of the Riemannian
manifolds $M_t=(M,g_t)$, where $g_t$ is the metric induced on
  $M$ by the map $\varphi_t$.\\
We say that a solution $\varphi_t$ is the flow by the gradient of the functional
  $\DG$ of the initial submanifold $\psi$.

By means of a slight extension of Polden's Theorem in~\cite{huiskpold}
(see~\cite{emin1} for details), there exists for some positive time a
unique smooth evolution $\varphi_t$ of any initial smooth submanifold
$M$. Our aim now is to show that under suitable hypotheses, such a flow
actually remains smooth for every time.

\begin{teo} If the differentiation order $k$ is strictly larger than
  $\left[\frac{n}{2}\right]+2$, then the flow by the gradient of
  $\DG$ of any initial $n$--dimensional submanifold is smooth for
  every positive time.\\
  Moreover, as $t\to+\infty$, the evolving submanifolds $\varphi_t$
  sub--converges (up to reparametrization and translation) to a
  smooth critical point of the functional $\DG$.
\end{teo}

Since the flow $\varphi_t$ is variational, the value of the functional
$\DG$ is monotone non increasing in time, hence it is bounded by its
value on the initial submanifold. This implies that, for all the
evolving submanifolds,
\begin{equation*}
{\mathrm {Vol}}(M_t)+\int\limits_{M}
\vert A^k\vert^2\,d\mu_t\leq C\,.
\end{equation*}
Hence, by means of Proposition~\ref{akteo} we get
\begin{equation*}
{\mathrm {Vol}}(M_t)+\int\limits_{M}
\vert\BBB\vert^{2k-4}\,d\mu_t\leq C
\end{equation*}
for a constant $C$ independent of time.\\
Since when $k>\left[\frac{n}{2}\right]+2$ we have $2k-4\geq n+1$ we
conclude that
\begin{equation}\label{npiuunobound}
{\mathrm {Vol}}(M_t)+\Vert\HHH\Vert_{L^{n+1}(\mu_t)}\leq C
\end{equation}
uniformly in time, for a constant $C$ depending only on the initial
submanifold.

By the results of Sections~5 and~6 in~\cite{mant5}, the
above a priori bound implies the following time--independent
interpolation inequalities for functions and tensors on $M$.

\begin{prop}\label{interp} As long as the flow $\varphi_t$ exists, 
for every smooth covariant tensor
$T=T_{i_1 \dots i_l}$ we have the inequalities
\begin{equation*}
    \Vert\nabla^j T\Vert_{L^{p}{(\mu_t)}}\leq\,C\,\Vert 
T\Vert_{W^{s,q}{(\mu_t)}}^{a}\Vert
      T\Vert_{L^r{(\mu_t)}}^{1-a}\,,
\end{equation*}
for all $j\in[0,s]$, $p, q, r\in[1,+\infty)$ and $a\in[j/s,1]$
with the compatibility condition
$$
\frac{1}{p}=\frac{j}{n}+a\left(\frac{1}{q}-\frac{s}{n}\right)+\frac{1-a}{r}\,.
$$
If such condition gives a negative value for $p$, the inequality holds
for every $p\in[1,+\infty)$ on the left side.\\
The constant $C$ depends on the dimensions, the  
orders of differentiation, the exponents of the involved norms and the
value of $\DG$ on the initial submanifold, 
but it is independent of time.
\end{prop}

Another consequence of inequality~\eqref{npiuunobound} 
is an uniform lower bound on the {\em Volume}
of $M_t$ (see the end of Section~5 in~\cite{mant5}), thus
\begin{equation*}
0 < c\leq {\mathrm {Vol}}(M_t)\leq C< +\infty\,,
\end{equation*}
again with a couple of constants $c$ and $C$ independent of time.

These a priori estimates allow us to forget the ``geometry'' of the
evolving submanifolds which, as it is well known, influences the
Sobolev constants (hence also the ones involved in the interpolation
inequalities). Thus, we can proceed with the estimates on the relevant
quantities that are the $L^2$ norms of the second fundamental form and its
derivatives.

From this point on, the rest of the proof follows step by step
Sections~7 and 8 in~\cite{mant5}, so we only give a sketch, referring 
to such paper for the details.

Differentiating in time the $L^2$ integrals
$\int_{M}\vert\nabla^{s}\BBB\vert^2\,d\mu_t$ during the flow, 
after some computations, one gets
$$
\frac{d\,}{dt}\int_M\vert\nabla^s\BBB\vert^2\,d\mu_t
\leq-3\int_M \vert\nabla^{s+k-1}\BBB\vert^2\,d\mu_t + 
\int_M \vert\qol_{2(s+k)}(\BBB)\vert + 
\vert\qol_{2(s+2)}(\BBB)\vert\,d\mu_t
$$
where $\qol_{2(s+k)}(\BBB)$ is a polynomial term not containing 
derivatives of $\BBB$ of order higher than $s+k-2$, such that any of
its monomials $Q_j(\BBB)$ has the following structure,
$$
Q_j(\BBB)= \nabla^{i_1}\BBB * \dots * \nabla^{i_N}\BBB\,,
$$
where the symbol $*$ means contraction with the scalar
product of $\R^{n+m}$, and $2(s+k)=\sum_{l=1}^N (i_l+1)$; the term 
$\qol_{2(s+2)}(\BBB)$ is analogous (see~\cite{mant5} for
more details).\\
Then, we can estimate them as follows
$$
\vert Q_j\vert  \leq \prod_{i=0}^{s+k-2}\vert\nabla^i \BBB\vert^{\alpha_{ji}} \qquad
\text{ with } \quad \sum_{i=0}^{s+k-2}\alpha_{ji}(i+1) = 2(s+k)\,.
$$
Then,
\begin{align*}
\int_{M} \vert Q_j\vert\,d\mu_t\leq \,&\,
\int_{M}\prod_{i=0}^{s+k-2}\vert\nabla^i
\BBB\vert^{\alpha_{ji}}\,d\mu_t\\
\leq\,&\,\prod_{i=0}^{s+k-2}\left(\int_{M}\vert\nabla^i
    \BBB\vert^{\alpha_{ji}\gamma_{i}}\,d\mu_t\right)^{\frac{1}{\gamma_i}}\\
=\,&\,\prod_{i=0}^{s+k-2}\Vert\nabla^i
\BBB\Vert_{L^{\alpha_{ji}\gamma_i}(\mu_t)}^{\alpha_{ij}}
\end{align*}
where the $\gamma_i$ are arbitrary positive values such that $\sum
1/{\gamma_i}=1$.\\
Interpolating now every factor of this product between
${\Vert\BBB\Vert}_{W^{s+k-1,2}(\mu_t)}$ and
${\Vert\BBB\Vert}_{L^{n+1}(\mu_t)}$ to some powers, by means of
Proposition~\ref{interp}, it is possible to conclude
$$
\int_{M} \vert Q_j\vert\,d\mu_t\leq\varepsilon_j\int_{M}\vert\nabla^{s+k-1}
\BBB\vert^2\,d\mu_t+C
$$
for arbitrarily small constants $\varepsilon_j>0$. Repeating this argument for all
the monomials $Q_j$ and choosing suitable $\varepsilon_j$ whose sum is
less than one,
$$
\frac{d\,}{dt}\int_{M}\vert\nabla^s\BBB\vert^2\,d\mu_t\,\leq
\,-2\int_{M}\vert\nabla^{s+k-1}\BBB\vert^2\,d\mu_t+C+
\int_M\vert\qol_{2(s+2)}(\BBB)\vert\,d\mu_t
$$
with a constant $C$ independent of time.\\
Dealing analogously with the other polynomial term
$\qol_{2(s+2)}(\BBB)$, we finally get
$$
\frac{d\,}{dt}\int_{M}\vert\nabla^s\BBB\vert^2\,d\mu_t\,\leq
\,-\int_{M}\vert\nabla^{s+k-1}\BBB\vert^2\,d\mu_t+C
\leq
\,-C\int_{M}\vert\nabla^{s}\BBB\vert^2\,d\mu_t+C\,,
$$
where in the second passage we used Poincar\`e inequality, which also
follows by Proposition~\ref{interp}.\\
Then, a simple ODE argument gives
$$
\int_{M}\vert\nabla^s \BBB\vert^2\,d\mu_t\leq C_s
$$
for every $s\in\NN$, with some constants $C_s$ dependent only on the
initial submanifold.\\
Again, via Proposition~\ref{interp}, one can now pass from this
family of Sobolev estimates to time--independent pointwise bounds
on all the covariant derivatives of $\BBB$.

Once we got these latter, the smoothness 
for every positive time and the sub--convergence follow by standard
arguments about geometric flows.

We underline the two key points where the properties of the distance
function play a role. First, when the order $k$ is
large enough, the estimate $\vert A^k\vert\geq
C_k\vert\BBB\vert^{2k-4}$ implies the a priori
estimates~\eqref{npiuunobound} leading to
the geometry--independent interpolation inequalities of
Proposition~\ref{interp}.\\
Secondly, the structure of $\vert A^k\vert^2$, 
(in particular the leading term, once expressed in terms
of the second fundamental form) produces a nice first variation
(computed in~\cite{ambman1}) giving rise to a well behaved parabolic
problem.

\bigskip

We conclude the paper mentioning the subsequent open problem suggested
by Ennio~De~Giorgi in the same paper~\cite{degio5} (Osservazione~2 and
Congettura~3, Pag.~267).

If $k>[n/2]+2$ it is easy to see, by the same proof, 
that for every $\varepsilon>0$ all the flows $\varphi^\varepsilon_t$,
associated to the functionals
$$
\DG^\varepsilon = \int_M
1+\varepsilon\vert\nabla^{k}\eta^M\vert^2\,d\mu\,,
$$
of a common initial $n$--dimensional submanifold, 
are smooth for every positive time.\\
Then, a natural question is the following:

\smallskip

{\em When the parameter $\varepsilon$ goes to zero, 
  the flows $\varphi^\varepsilon_t$ converge, in some sense, 
  to the flow associated to the limit functional which is simply the
  ${\mathrm {Volume}}$ functional?\\
  That is, do they converge to the {\em motion by mean curvature} of
  the initial submanifold?}

\bibliographystyle{amsplain}
\bibliography{dgcong}

\end{document}